\input amstex
\documentstyle{amsppt}

\def\stydate{June 26, 2000}

\chardef\tempcat\catcode`\@ \ifx\undefined\amstexloaded\input amstex
\else\catcode`\@\tempcat\fi \expandafter\ifx\csname
amsppt.sty\endcsname\relax\input amsppt.sty \fi
\let\tempcat\undefined

\immediate\write16{This is LABEL.DEF by A.Degtyarev <\stydate>}
\expandafter\ifx\csname label.def\endcsname\relax\else
  \message{[already loaded]}\endinput\fi
\expandafter\edef\csname label.def\endcsname{%
  \catcode`\noexpand\@\the\catcode`\@\edef\noexpand\styname{LABEL.DEF}%
  \def\expandafter\noexpand\csname label.def\endcsname{\stydate}%
    \toks0{}\toks2{}}
\catcode`\@11
\def\labelmesg@ {LABEL.DEF: }
{\edef\temp{\the\everyjob\W@{\labelmesg@<\stydate>}}
\global\everyjob\expandafter{\temp}}

\def\@car#1#2\@nil{#1}
\def\@cdr#1#2\@nil{#2}
\def\eat@bs{\expandafter\eat@\string}
\def\eat@ii#1#2{}
\def\eat@iii#1#2#3{}
\def\eat@iv#1#2#3#4{}
\def\@DO#1#2\@{\expandafter#1\csname\eat@bs#2\endcsname}
\def\@N#1\@{\csname\eat@bs#1\endcsname}
\def\@Nx{\@DO\noexpand}
\def\@Name#1\@{\if\@undefined#1\@\else\@N#1\@\fi}
\def\@Ndef{\@DO\def}
\def\@Ngdef{\global\@Ndef}
\def\@Nedef{\@DO\edef}
\def\@Nxdef{\global\@Nedef}
\def\@Nlet{\@DO\let}
\def\@undefined#1\@{\@DO\ifx#1\@\relax\true@\else\false@\fi}
\def\@@addto#1#2{{\toks@\expandafter{#1#2}\xdef#1{\the\toks@}}}
\def\@@addparm#1#2{{\toks@\expandafter{#1{##1}#2}%
    \edef#1{\gdef\noexpand#1####1{\the\toks@}}#1}}
\def\make@letter{\edef\t@mpcat{\catcode`\@\the\catcode`\@}\catcode`\@11 }
\def\donext@{\expandafter\egroup\next@}
\def\x@notempty#1{\expandafter\notempty\expandafter{#1}}
\def\lc@def#1#2{\edef#1{#2}%
    \lowercase\expandafter{\expandafter\edef\expandafter#1\expandafter{#1}}}
\newif\iffound@
\def\find@#1\in#2{\found@false
    \DNii@{\ifx\next\@nil\let\next\eat@\else\let\next\nextiv@\fi\next}%
    \edef\nextiii@{#1}\def\nextiv@##1,{%
    \edef\next{##1}\ifx\nextiii@\next\found@true\fi\FN@\nextii@}%
    \expandafter\nextiv@#2,\@nil}
{\let\head\relax\let\specialhead\relax\let\subhead\relax
\let\subsubhead\relax\let\proclaim\relax
\gdef\let@relax{\let\head\relax\let\specialhead\relax\let\subhead\relax
    \let\subsubhead\relax\let\proclaim\relax}}
\newskip\@savsk
\let\@ignorespaces\ignorespaces
\def\@ignorespacesp{\ifhmode
  \ifdim\lastskip>\z@\else\penalty\@M\hskip-1sp%
        \penalty\@M\hskip1sp \fi\fi\@ignorespaces}
\def\ignorespaces{\protect\@ignorespacesp}
\def\@bsphack{\relax\ifmmode\else\@savsk\lastskip
  \ifhmode\edef\@sf{\spacefactor\the\spacefactor}\fi\fi}
\def\@esphack{\relax
  \ifx\penalty@\penalty\else\penalty\@M\fi   
  \ifmmode\else\ifhmode\@sf{}\ifdim\@savsk>\z@\@ignorespacesp\fi\fi\fi}
\let\@frills@\identity@
\let\@txtopt@\identyty@
\newif\if@star
\newif\if@write\@writetrue
\def\@numopt@{\if@star\expandafter\eat@\fi}
\def\checkstar@#1{\DN@{\@writetrue
  \ifx\next*\DN@####1{\@startrue\checkstar@@{#1}}%
      \else\DN@{\@starfalse#1}\fi\next@}\FN@\next@}
\def\checkstar@@#1{\DN@{%
  \ifx\next*\DN@####1{\@writefalse#1}%
      \else\DN@{\@writetrue#1}\fi\next@}\FN@\next@}
\def\checkfrills@#1{\DN@{%
  \ifx\next\nofrills\DN@####1{#1}\def\@frills@####1{####1\nofrills}%
      \else\DN@{#1}\let\@frills@\identity@\fi\next@}\FN@\next@}
\def\checkbrack@#1{\DN@{%
    \ifx\next[\DN@[####1]{\def\@txtopt@########1{####1}#1}%
    \else\DN@{\let\@txtopt@\identity@#1}\fi\next@}\FN@\next@}
\def\check@therstyle#1#2{\bgroup\DN@{#1}\ifx\@txtopt@\identity@\else
        \DNii@##1\@therstyle{}\def\@therstyle{\DN@{#2}\nextii@}%
    \expandafter\expandafter\expandafter\nextii@\@txtopt@\@therstyle.\@therstyle
    \fi\donext@}

\newread\@inputcheck
\def\@input#1{\openin\@inputcheck #1 \ifeof\@inputcheck \W@
  {No file `#1'.}\else\closein\@inputcheck \relax\input #1 \fi}

\def\loadstyle#1{\edef\next{#1}%
    \DN@##1.##2\@nil{\if\notempty{##2}\else\def\next{##1.sty}\fi}%
    \expandafter\next@\next.\@nil\lc@def\next@\next
    \expandafter\ifx\csname\next@\endcsname\relax\input\next\fi}

\let\pagebody@\pagebody
\let\pagetop@\empty
\let\pagebot@\empty
\let\@Xend\empty
\def\pagebody{\pagetop@\pagebody@\pagebot@\@Xend}
\let\@Xclose\empty

\newwrite\@Xmain
\newwrite\@Xsub
\def\W@X{\write\@Xout}
\def\make@Xmain{\global\let\@Xout\@Xmain\global\let\end\endmain@
  \xdef\@Xname{\jobname}\xdef\@inputname{\jobname}}
\begingroup
\catcode`\(\the\catcode`\{\catcode`\{12
\catcode`\)\the\catcode`\}\catcode`\}12
\gdef\W@count#1((\lc@def\@tempa(#1)%
    \def\\##1(\W@X(\global##1\the##1))%
    \edef\@tempa(\W@X(%
        \string\expandafter\gdef\string\csname\space\@tempa\string\endcsname{)%
        \\\pageno\\\footmarkcount@\@Xclose\W@X(}))\expandafter)\@tempa)
\endgroup
\def\readaux{\bgroup\checkbrack@\readaux@}
\let\begin\readaux
\def\readaux@{%
    \W@{>>> \labelmesg@ Run this file twice to get x-references right}%
    \global\everypar{}%
    {\def\\{\global\let}%
        \def\/##1##2{\gdef##1{\wrn@command##1##2}}%
        \disablepreambule@cs}%
    \make@Xmain{\make@letter\setboxz@h{\@input{\@txtopt@{\@Xname.aux}}%
            \lc@def\@tempa\jobname\@Name\open@\@tempa\@}}%
  \immediate\openout\@Xout\@Xname.aux%
    \immediate\W@X{\relax}\egroup}
\everypar{\global\everypar{}\readaux} {\toks@\expandafter{\topmatter}
\global\edef\topmatter{\noexpand\readaux\the\toks@}}
\let\@@end@@\end

\def\@Xclose@{{\def\@Xend{\ifnum\insertpenalties=\z@
        \W@count{close@\@Xname}\closeout\@Xout\fi}%
    \vfill\supereject}}
\def\endmain@{\@Xclose@
    \W@{>>> \labelmesg@ Run this file twice to get x-references right}%
    \@@end@@}
\def\disablepreambule@cs{\\\disablepreambule@cs\relax}

\def\include#1{\bgroup
  \ifx\@Xout\@Xsub\DN@{\errmessage
        {\labelmesg@ Only one level of \string\include\space is supported}}%
    \else\edef\@tempb{#1}\clearpage
      \DN@##1 {\if\notempty{##1}\edef\@tempb{##1}\DN@####1\eat@ {}\fi\next@}%
    \DNii@##1.{\edef\@tempa{##1}\DN@####1\eat@.{}\next@}%
        \expandafter\next@\@tempb\eat@{} \eat@{} %
    \expandafter\nextii@\@tempb.\eat@.%
        \relaxnext@
      \if\x@notempty\@tempa
          \edef\nextii@{\write\@Xmain{%
            \noexpand\string\noexpand\@input{\@tempa.aux}}}\nextii@
        \ifx\undefined\@includelist\found@true\else
                    \find@\@tempa\in\@includelist\fi
            \iffound@\ifx\undefined\@noincllist\found@false\else
                    \find@\@tempb\in\@noincllist\fi\else\found@true\fi
            \iffound@\lc@def\@tempa\@tempa
                \if\@undefined\close@\@tempa\@\else\edef\next@{\@Nx\close@\@tempa\@}\fi
            \else\xdef\@Xname{\@tempa}\xdef\@inputname{\@tempb}%
                \W@count{open@\@Xname}\global\let\@Xout\@Xsub
            \openout\@Xout\@tempa.aux \W@X{\relax}%
            \DN@{\let\end\endinput\@input\@inputname
                    \@Xclose@\make@Xmain}\fi\fi\fi
  \donext@}
\def\includeonly#1{\edef\@includelist{#1}}
\def\noinclude#1{\edef\@noincllist{#1}}

\def\arabicnum#1{\number#1}

\def\Romannum#1{\expandafter\uppercase\expandafter{\romannumeral#1}}
\def\alphnum#1{\ifcase#1\or a\or b\or c\or d\else\@ialph{#1}\fi}
\def\@ialph#1{\ifcase#1\or \or \or \or \or e\or f\or g\or h\or i\or j\or
    k\or l\or m\or n\or o\or p\or q\or r\or s\or t\or u\or v\or w\or x\or y\or
    z\else\fi}
\def\Alphnum#1{\ifcase#1\or A\or B\or C\or D\else\@Ialph{#1}\fi}
\def\@Ialph#1{\ifcase#1\or \or \or \or \or E\or F\or G\or H\or I\or J\or
    K\or L\or M\or N\or O\or P\or Q\or R\or S\or T\or U\or V\or W\or X\or Y\or
    Z\else\fi}

\def\ST@P{step}
\def\ST@LE{style}
\def\N@M{no}
\def\F@NT{font@}
\outer\def\newcounter{\checkbrack@{\expandafter\newcounter@\@txtopt@{{}}}}
{\let\newcount\relax
\gdef\newcounter@#1#2#3{{%
    \toks@@\expandafter{\csname\eat@bs#2\N@M\endcsname}%
    \DN@{\alloc@0\count\countdef\insc@unt}%
    \ifx\@txtopt@\identity@\expandafter\next@\the\toks@@
        \else\if\notempty{#1}\global\@Nlet#2\N@M\@#1\fi\fi
    \@Nxdef\the\eat@bs#2\@{\if\@undefined\the\eat@bs#3\@\else
            \@Nx\the\eat@bs#3\@.\fi\noexpand\arabicnum\the\toks@@}%
  \@Nxdef#2\ST@P\@{}%
  \if\@undefined#3\ST@P\@\else
    \edef\next@{\noexpand\@@addto\@Nx#3\ST@P\@{%
             \global\@Nx#2\N@M\@\z@\@Nx#2\ST@P\@}}\next@\fi
    \expandafter\@@addto\expandafter\@Xclose\expandafter
        {\expandafter\\\the\toks@@}}}}
\outer\def\copycounter#1#2{%
    \@Nxdef#1\N@M\@{\@Nx#2\N@M\@}%
    \@Nxdef#1\ST@P\@{\@Nx#2\ST@P\@}%
    \@Nxdef\the\eat@bs#1\@{\@Nx\the\eat@bs#2\@}}
\outer\def\everystep{\checkstar@\everystep@}
\def\everystep@#1{\if@star\let\next@\gdef\else\let\next@\@@addto\fi
    \@DO\next@#1\ST@P\@}
\def\counterstyle#1{\@Ngdef\the\eat@bs#1\@}
\def\advancecounter#1#2{\@N#1\ST@P\@\global\advance\@N#1\N@M\@#2}
\def\setcounter#1#2{\@N#1\ST@P\@\global\@N#1\N@M\@#2}
\def\counter#1{\refstepcounter#1\printcounter#1}
\def\printcounter#1{\@N\the\eat@bs#1\@}
\def\refcounter#1{\xdef\@lastmark{\printcounter#1}}
\def\stepcounter#1{\advancecounter#1\@ne}
\def\refstepcounter#1{\stepcounter#1\refcounter#1}
\def\savecounter#1{\@Nedef#1@sav\@{\global\@N#1\N@M\@\the\@N#1\N@M\@}}
\def\restorecounter#1{\@Name#1@sav\@}

\def\warning#1#2{\W@{Warning: #1 on input line #2}}
\def\warning@#1{\warning{#1}{\the\inputlineno}}
\def\wrn@@Protect#1#2{\warning@{\string\Protect\string#1\space ignored}}
\def\wrn@@label#1#2{\warning{label `#1' multiply defined}{#2}}
\def\wrn@@ref#1#2{\warning@{label `#1' undefined}}
\def\wrn@@cite#1#2{\warning@{citation `#1' undefined}}
\def\wrn@@command#1#2{\warning@{Preamble command \string#1\space ignored}#2}
\def\wrn@@option#1#2{\warning@{Option \string#1\string#2\space is not supported}}
\def\wrn@@reference#1#2{\W@{Reference `#1' on input line \the\inputlineno}}
\def\wrn@@citation#1#2{\W@{Citation `#1' on input line \the\inputlineno}}
\let\wrn@reference\eat@ii
\let\wrn@citation\eat@ii
\def\nowarning#1{\if\@undefined\wrn@\eat@bs#1\@\wrn@option\nowarning#1\else
        \@Nlet\wrn@\eat@bs#1\@\eat@ii\fi}
\def\printwarning#1{\if\@undefined\wrn@@\eat@bs#1\@\wrn@option\printwarning#1\else
        \@Nlet\wrn@\eat@bs#1\expandafter\@\csname wrn@@\eat@bs#1\endcsname\fi}
\printwarning\Protect \printwarning\label \printwarning\ref
\printwarning\cite \printwarning\command \printwarning\option

\def\bftext{\ifmmode\fam\bffam\else\bf\fi}
\let\@lastmark\empty
\let\@lastlabel\empty
\def\lastmark{\@lastmark}
\let\lastlabel\empty
\let\everylabel\relax
\let\everylabel@\eat@
\let\everyref\relax
\def\newlabel{\bgroup\everylabel\newlabel@}
\def\newlabel@#1#2#3{\if\@undefined\r@-#1\@\else\wrn@label{#1}{#3}\fi
  {\let\protect\noexpand\@Nxdef\r@-#1\@{#2}}\egroup}
\def\w@ref{\bgroup\everyref\w@@ref}
\def\w@@ref#1#2#3#4{%
    \if\@undefined\r@-#1\@{\bftext??}#2{#1}{}\else
  \@DO{\expandafter\expandafter#3}\r@-#1\@\@nil\fi#4{#1}{}\egroup}
\def\@@@xref#1{\w@ref{#1}\wrn@ref\@car\wrn@reference}
\def\@xref#1{\rom{\@@@xref{#1}}}
\let\xref\@xref
\def\pageref#1{\w@ref{#1}\wrn@ref\@cdr\wrn@reference}
\def\thepage{\ifnum\pageno<\z@\romannumeral-\pageno\else\number\pageno\fi}
\def\label@{\@bsphack\bgroup\everylabel\label@@}
\def\label@@#1#2{\everylabel@{{#1}{#2}}\let\thepage\relax
  \def\protect{\noexpand\noexpand\noexpand}%
  \edef\@tempa{\edef\noexpand\@lastlabel{#1}%
    \W@X{\string\newlabel{#2}{{\@lastmark}{\thepage}}{\the\inputlineno}}}%
  \expandafter\egroup\@tempa\@esphack}
\def\label#1{\label@{#1}{#1}}
\def\fn@P@{\relaxnext@
    \DN@{\ifx[\next\DN@[####1]{}\else
        \ifx"\next\DN@"####1"{}\else\DN@{}\fi\fi\next@}%
    \FN@\next@}
\def\eat@fn#1{\ifx#1[\expandafter\eat@br\else
  \ifx#1"\expandafter\expandafter\expandafter\eat@qu\fi\fi}
\def\eat@br#1]#2{}
\def\eat@qu#1"#2{}
{\catcode`\~\active\lccode`\~`\@
\lowercase{\global\let\@@P@~\gdef~{\protect\@@P@}}}
\def\Protect@@#1{\def#1{\protect#1}}
\def\disable@special{\let\W@X@\eat@iii\let\label\eat@
    \def\footnotemark{\protect\fn@P@}%
  \let\footnotetext\eat@fn\let\footnote\eat@fn
    \let\refcounter\eat@\let\savecounter\eat@\let\restorecounter\eat@
    \let\advancecounter\eat@ii\let\setcounter\eat@ii
  \let\ifvmode\iffalse\Protect@@\@@@xref\Protect@@\pageref\Protect@@\nofrills
    \Protect@@\\\Protect@@~}
\let\notoctext\identity@
\def\W@X@#1#2#3{\@bsphack{\disable@special\let\notoctext\eat@
    \def\chapter{\protect\chapter@toc}\let\thepage\relax
    \def\protect{\noexpand\noexpand\noexpand}#1%
  \edef\next@{\if\@undefined#2\@\else\write#2{#3}\fi}\expandafter}\next@
    \@esphack}
\def\writeauxline#1#2#3{\W@X@{}\@Xout{\string\@Xline{#1}{#2}{#3}{\thepage}}}
{\let\newwrite\relax
\gdef\@openin#1{\make@letter\@input{\jobname.#1}\t@mpcat}
\gdef\@openout#1{\global\expandafter\newwrite\csname tf@-#1\endcsname
   \immediate\openout\@N\tf@-#1\@\jobname.#1\relax}}
\def\@@openout#1{\@openout{#1}%
  \@@addto\readaux@{\immediate\closeout\@N\tf@-#1\@}}
\def\auxlinedef#1{\@Ndef\do@-#1\@}
\def\@Xline#1{\if\@undefined\do@-#1\@\expandafter\eat@iii\else
    \@DO\expandafter\do@-#1\@\fi}
\def\beginW@{\bgroup\def\do##1{\catcode`##112 }\dospecials\do\@\do\"
    \catcode`\{\@ne\catcode`\}\tw@\immediate\write\@N}
\def\endW@toc#1#2#3{{\string\tocline{#1}{#2\string\page{#3}}}\egroup}
\def\do@tocline#1{%
    \if\@undefined\tf@-#1\@\expandafter\eat@iii\else
        \beginW@\tf@-#1\@\expandafter\endW@toc\fi
} \auxlinedef{toc}{\do@tocline{toc}}

\let\protect\empty
\def\Protect#1{\if\@undefined#1@P@\@\PROTECT#1\else\wrn@Protect#1\empty\fi}
\def\PROTECT#1{\@Nlet#1@P@\@#1\edef#1{\noexpand\protect\@Nx#1@P@\@}}
\def\pdef#1{\edef#1{\noexpand\protect\@Nx#1@P@\@}\@Ndef#1@P@\@}

\Protect\operatorname \Protect\operatornamewithlimits \Protect\qopname@
\Protect\qopnamewl@ \Protect\text \Protect\topsmash \Protect\botsmash
\Protect\smash \Protect\widetilde \Protect\widehat \Protect\thetag
\Protect\therosteritem
\Protect\Cal \Protect\Bbb \Protect\bold \Protect\slanted \Protect\roman
\Protect\italic \Protect\boldkey \Protect\boldsymbol \Protect\frak
\Protect\goth \Protect\dots
\Protect\cong \Protect\lbrace \let\{\lbrace \Protect\rbrace \let\}\rbrace
\let\root@P@@\root \def\root@P@#1{\root@P@@#1\of}
\def\root#1\of{\protect\root@P@{#1}}

\def\frills{\ignorespaces\@txtopt@}
\def\frillsnotempty#1{\x@notempty{\@txtopt@{#1}}}
\def\numberline{\@numopt@}
\newif\if@theorem
\let\@therstyle\eat@
\def\@headtext@#1#2{{\disable@special\let\protect\noexpand
    \def\chapter{\protect\chapter@rh}%
    \edef\next@{\noexpand\@frills@\noexpand#1{#2}}\expandafter}\next@}
\let\AmSrighthead@\rightheadtext
\def\rightheadtext{\checkfrills@{\@headtext@\AmSrighthead@}}
\let\AmSlefthead@\leftheadtext
\def\leftheadtext{\checkfrills@{\@headtext@\AmSlefthead@}}
\def\@head@@#1#2#3#4#5{\@Name\pre\eat@bs#1\@\if@theorem\else
    \@frills@{\csname\expandafter\eat@iv\string#4\endcsname}\relax
        \ifx\protect\empty\@N#1\F@NT\@\fi\fi
    \@N#1\ST@LE\@{\counter#3}{#5}%
  \if@write\writeauxline{toc}{\eat@bs#1}{#2{\counter#3}{#5}}\fi
    \if@theorem\else\expandafter#4\fi
    \ifx#4\endhead\ifx\@txtopt@\identity@\else
        \headmark{\@N#1\ST@LE\@{\counter#3}{\frills\empty}}\fi\fi
    \@Name\post\eat@bs#1\@\ignorespaces}
\ifx\undefined\endhead\Invalid@\endhead\fi
\def\@head@#1{\checkstar@{\checkfrills@{\checkbrack@{\@head@@#1}}}}
\def\@thm@@#1#2#3{\@Name\pre\eat@bs#1\@
    \@frills@{\csname\expandafter\eat@iv\string#3\endcsname}
    {\@theoremtrue\check@therstyle{\@N#1\ST@LE\@}\frills
            {\counter#2}\@theoremfalse}%
    \@DO\envir@stack\end\eat@bs#1\@
    \@N#1\F@NT\@\@Name\post\eat@bs#1\@\ignorespaces}
\def\@thm@#1{\checkstar@{\checkfrills@{\checkbrack@{\@thm@@#1}}}}
\def\@capt@@#1#2#3#4#5\endcaption{\bgroup
    \edef\@tempb{\global\footmarkcount@\the\footmarkcount@
    \global\@N#2\N@M\@\the\@N#2\N@M\@}%
    \def\shortcaption##1{\global\def\sh@rtt@xt####1{##1}}\let\sh@rtt@xt\identity@
    \DN@{#4{\@tempb\@N#1\ST@LE\@{\counter#2}}}%
    \if\notempty{#5}\DNii@{\next@\@N#1\F@NT\@}\else\let\nextii@\next@\fi
    \nextii@#5\endcaption
  \if@write\writeauxline{#3}{\eat@bs#1}{{} \@N#1\ST@LE\@{\counter#2}%
    \if\notempty{#5}.\enspace\fi\sh@rtt@xt{#5}}\fi
  \global\let\sh@rtt@xt\undefined\egroup}
\def\@capt@#1{\checkstar@{\checkfrills@{\checkbrack@{\@capt@@#1}}}}
\let\captiontextfont@\empty

\ifx\undefined\subsubheadfont@\def\subsubheadfont@{\it}\fi
\ifx\undefined\proclaimfont\def\proclaimfont{\sl}\fi
\ifx\undefined\proclaimfont@\let\proclaimfont@\proclaimfont\fi
\def\proclaimfont{\proclaimfont@}
\ifx\undefined\definitionfont@\def\AmSdeffont@{\rm}
    \else\let\AmSdeffont@\definitionfont@\fi
\ifx\undefined\remarkfont@\def\remarkfont@{\rm}\fi

\def\newfont@def#1#2{\if\@undefined#1\F@NT\@
    \@Nxdef#1\F@NT\@{\@Nx.\expandafter\eat@iv\string#2\F@NT\@}\fi}
\def\newhead@#1#2#3#4{{%
    \gdef#1{\@therstyle\@therstyle\@head@{#1#2#3#4}}\newfont@def#1#4%
    \if\@undefined#1\ST@LE\@\@Ngdef#1\ST@LE\@{\headstyle}\fi
    \if\@undefined#2\@\gdef#2{\headtocstyle}\fi
  \@@addto\moretocdefs@{\\#1#1#4}}}
\outer\def\newhead#1{\checkbrack@{\expandafter\newhead@\expandafter
    #1\@txtopt@\headtocstyle}}
\outer\def\newtheorem#1#2#3#4{{%
    \gdef#2{\@thm@{#2#3#4}}\newfont@def#2#4%
    \@Nxdef\end\eat@bs#2\@{\noexpand\revert@envir
        \@Nx\end\eat@bs#2\@\noexpand#4}%
  \if\@undefined#2\ST@LE\@\@Ngdef#2\ST@LE\@{\proclaimstyle{#1}}\fi}}%
\outer\def\newcaption#1#2#3#4#5{{\let#2\relax
  \edef\@tempa{\gdef#2####1\@Nx\end\eat@bs#2\@}%
    \@tempa{\@capt@{#2#3{#4}#5}##1\endcaption}\newfont@def#2\endcaptiontext%
  \if\@undefined#2\ST@LE\@\@Ngdef#2\ST@LE\@{\captionstyle{#1}}\fi
  \@@addto\moretocdefs@{\\#2#2\endcaption}\newtoc{#4}}}
{
\outer\gdef\newtoc#1{{%
    \@DO\ifx\do@-#1\@\relax
    \global\auxlinedef{#1}{\do@tocline{#1}}{}%
    \@@addto\tocsections@{\make@toc{#1}{}}\fi}}}

\toks@\expandafter{\itembox@}
\toks@@{\bgroup\let\therosteritem\identity@\let\rm\empty
  \edef\next@{\edef\noexpand\@lastmark{\therosteritem@}}\donext@}
\edef\itembox@{\the\toks@@\the\toks@}
\def\firstitem@false{\let\iffirstitem@\iffalse
    \global\let\lastlabel\@lastlabel}

\let\rosteritemrefform\therosteritem
\let\rosteritemrefseparator\empty
\def\rosteritemref#1{\hbox{\rosteritemrefform{\@@@xref{#1}}}}
\def\local#1{\label@\@lastlabel{\lastlabel-i#1}}
\def\loccit#1{\rosteritemref{\lastlabel-i#1}}
\def\xRef@P@{\gdef\lastlabel}
\def\xRef#1{\@xref{#1}\protect\xRef@P@{#1}}

\def\iref@P@{\gdef\lastref}
\def\itemref#1#2{\rosteritemref{#1-i#2}\protect\iref@P@{#1}}
\def\iref#1{\@xref{#1}\rosteritemrefseparator\itemref{#1}}
\def\ditto#1{\rosteritemref{\lastref-i#1}}

\def\eqref#1{\thetag{\@@@xref{#1}}}
\def\tagform@#1{\ifmmode\hbox{\rm\else\rom{\fi
        (\ignorespaces#1\unskip)\iftrue}\else}\fi}

\let\AmSfnote@\makefootnote@
\def\makefootnote@#1{\bgroup\let\footmarkform@\identity@
  \edef\next@{\edef\noexpand\@lastmark{#1}}\donext@\AmSfnote@{#1}}

\def\clearpage{\ifnum\insertpenalties>0\line{}\fi\vfill\supereject}

\def\proof{\checkfrills@{\checkbrack@{%
    \check@therstyle{\@frills@{\demo}{\frills{Proof}}{}}
        {\frills{}\envir@stack\endremark\envir@stack\enddemo}%
  \envir@stack\endproof\ignorespaces}}}
\def\endproof{\nofrillscheck{\frills@{\qed}\revert@envir\endproof\enddemo}}

\let\AmSref\ref
\let\AmSrefstyle\refstyle
\let\plaincite\cite
\def\citei@#1,{\citeii@#1\eat@,}
\def\citeii@#1\eat@{\w@ref{#1}\wrn@cite\@car\wrn@citation}
\def\mcite@#1;{\plaincite{\citei@#1\eat@,\unskip}\mcite@i}
\def\mcite@i#1;{\DN@{#1}\ifx\next@\endmcite@
  \else, \plaincite{\citei@#1\eat@,\unskip}\expandafter\mcite@i\fi}
\def\endmcite@{\endmcite@}
\def\cite#1{\mcite@#1;\endmcite@;}
\PROTECT\cite
\def\refstyle#1{\AmSrefstyle{#1}\uppercase{%
    \ifx#1A\relax \def\@ref@##1{\AmSref\xdef\@lastmark{##1}\key##1}%
    \else\ifx#1C\relax \def\@ref@{\AmSref\no\counter\refno}%
        \else\def\@ref@{\AmSref}\fi\fi}}
\refstyle A
\newcounter\refno\null
\newif\ifRefs
\gdef\Refs{\checkstar@{\checkbrack@{\csname AmSRefs\endcsname
  \nofrills{\frills{References}%
  \if@write\writeauxline{toc}{vartocline}{\frills{References}}\fi}%
  \def\ref{\@ref@}\Refstrue\ignorespaces}}}
\let\ref\xref

\newif\iftoc
\pdef\tocbreak{\iftoc\hfil\break\fi}
\def\tocsections@{\make@toc{toc}{}}
\let\moretocdefs@\empty
\def\newtocline@#1#2#3{%
  \edef#1{\def\@Nx#2line\@####1{\@Nx.\expandafter\eat@iv
        \string#3\@####1\noexpand#3}}%
  \@Nedef\no\eat@bs#1\@{\let\@Nx#2line\@\noexpand\eat@}%
    \@N\no\eat@bs#1\@}
\def\MakeToc#1{\@@openout{#1}}
\def\newtocline#1#2#3{\Err@{\Invalid@@\string\newtocline}}
\def\make@toc#1#2{\penaltyandskip@{-200}\aboveheadskip
    \if\notempty{#2}
        \centerline{\headfont@\ignorespaces#2\unskip}\nobreak
    \vskip\belowheadskip \fi
    \@openin{#1}\relax
    \vskip\z@}
\def\contents{\readaux\checkfrills@{\checkbrack@{\@contents@}}}
\def\@contents@{\toc@{\frills{Contents}}\envir@stack\endcontents%
    \def\nopagenumbers{\let\page\eat@}\let\newtocline\newtocline@\toctrue
  \def\tocline##1{\csname##1line\endcsname}%
  \edef\caption##1\endcaption{\expandafter\noexpand
    \csname head\endcsname##1\noexpand\endhead}%
    \ifmonograph@\def\vartoclineline{\Chapterline}%
        \else\def\vartoclineline##1{\sectionline{{} ##1}}\fi
  \let\\\newtocline@\moretocdefs@
    \ifx\@frills@\identity@\def\\##1##2##3{##1}\moretocdefs@
        \else\let\tocsections@\relax\fi
    \def\\{\unskip\space\ignorespaces}\let\maketoc\make@toc}
\def\endcontents{\tocsections@\vskip-\lastskip\revert@envir\endcontents
    \endtoc}

\if\@undefined\selectf@nt\@\let\selectf@nt\identity@\fi
\def\Err@math#1{\Err@{Use \string#1\space only in text}}
\def\textonlyfont@#1#2{%
    \def#1{\RIfM@\Err@math#1\else\edef\f@ntsh@pe{\string#1}\selectf@nt#2\fi}%
    \PROTECT#1}
\tenpoint

\def\newshapeswitch#1#2{\gdef#1{\selectsh@pe#1#2}\PROTECT#1}
\def\shapeswitch#1#2#3{\@Ngdef#1\string#2\@{#3}}
\shapeswitch\rm\bf\bf \shapeswitch\rm\tt\tt \shapeswitch\rm\smc\smc
\newshapeswitch\em\it
\shapeswitch\em\it\rm \shapeswitch\em\sl\rm
\def\selectsh@pe#1#2{\relax\if\@undefined#1\f@ntsh@pe\@#2\else
    \@N#1\f@ntsh@pe\@\fi}

\def\@itcorr@{\leavevmode
    \edef\prevskip@{\ifdim\lastskip=\z@ \else\hskip\the\lastskip\relax\fi}\unskip
    \edef\prevpenalty@{\ifnum\lastpenalty=\z@ \else
        \penalty\the\lastpenalty\relax\fi}\unpenalty
    \/\prevpenalty@\prevskip@}
\def\rom@P@#1{\@itcorr@{\selectsh@pe\rm\rm#1}}
\def\rom{\protect\rom@P@}
\def\Rom@P@#1{\@itcorr@{\rm#1}}
\def\Rom{\protect\Rom@P@}
{\catcode`\-11
\newcount\cnt@-idx \global\cnt@-idx=10000
\newcount\cnt@-glo \global\cnt@-glo=10000
\gdef\writeindex#1{\W@X@{\global\advance\cnt@-idx\@ne}\tf@-idx
 {\string\indexentry{#1}{\the\cnt@-idx}{\thepage}}}
\gdef\writeglossary#1{\W@X@{\global\advance\cnt@-glo\@ne}\tf@-glo
 {\string\glossaryentry{#1}{\the\cnt@-glo}{\thepage}}}}
\def\emph#1{\@itcorr@\bgroup\em\ignorespaces#1\unskip\egroup
  \DN@{\DN@{}\ifx\next.\else\ifx\next,\else\DN@{\/}\fi\fi\next@}\FN@\next@}
\def\makequoteactive{\catcode`\"\active}
{\makequoteactive\gdef"{\FN@\quote@}
\gdef\quote@{\ifx"\next\DN@"##1""{\quoteii{##1}}\else\DN@##1"{\quotei{##1}}\fi\next@}}
\let\quotei\eat@
\let\quoteii\eat@
\def\MakeIndex{\@openout{idx}}
\def\MakeGlossary{\@openout{glo}}

\def\endofpar#1{\ifmmode\ifinner\endofpar@{#1}\else\eqno{#1}\fi
    \else\leavevmode\endofpar@{#1}\fi}
\def\endofpar@#1{\unskip\penalty\z@\null\hfil\hbox{#1}\hfilneg\penalty\@M}

\newdimen\normalparindent\normalparindent\parindent
\def\firstparindent#1{\everypar\expandafter{\the\everypar
  \global\parindent\normalparindent\global\everypar{}}\parindent#1\relax}

\@@addto\disablepreambule@cs{%
    \\\readaux\relax
    \\\begin\relax
    \\\readaux@\relax
    \\\@openout\eat@
    \\\@@openout\eat@
    \/\Monograph\empty
    \/\MakeIndex\empty
    \/\MakeGlossary\empty
    \/\MakeToc\eat@
}

\csname label.def\endcsname


\def\punct#1#2{\if\notempty{#2}#1\fi}
\def\sppunct{\punct{.\enspace}}
\def\varpunct#1#2{\if\frillsnotempty{#2}#1\fi}

\def\headstyle#1#2{\numberline{#1\sppunct{#2}}\ignorespaces#2\unskip}
\def\headtocstyle#1#2{\numberline{#1\punct.{#2}}\space #2}

\def\specialtocstyle#1#2{#2}
\newcounter\section\null
\newcounter\subsection\section
\newcounter\subsubsection\subsection
\newhead\specialsection[\specialtocstyle]\null\endspecialhead
\newhead\section\section\endhead
\newhead\subsection\subsection\endsubhead
\newhead\subsubsection\subsubsection\endsubsubhead
\def\firstappendix{\global\sectionno0 %
  \counterstyle\section{\Alphnum\sectionno}%
    \global\let\firstappendix\empty}

\def\appendixtocstyle#1#2{\space\numberline{Appendix #1\sppunct{#2}}#2}
\newhead\appendix[\appendixtocstyle]\section\endhead

\let\endAmSdef\enddefinition
\def\proclaimstyle#1#2{\numberline{#2\varpunct{.\enspace}{#1}}\frills{#1}}
\copycounter\thm\subsubsection
\theorem\thm\endproclaim
\proposition\thm\endproclaim
\lemma\thm\endproclaim
\corollary\thm\endproclaim
\definition\thm\endAmSdef
\example\thm\endAmSdef

\def\captionstyle#1#2{\frills{#1}\numberline{\varpunct{ }{#1}#2}}
\newcounter\figure\null
\newcounter\table\null
\newcaption{Figure}\figure\figure{lof}\botcaption
\newcaption{Table}\table\table{lot}\topcaption

\copycounter\equation\subsubsection


\expandafter\ifx\csname label.def\endcsname\relax\input label.def \fi
\def\stydate{June 26, 2000}
\def\styname{DEBUG.DEF}
\immediate\write16{This is \styname\space by A.Degtyarev <\stydate>}
\expandafter\ifx\csname debug.def\endcsname\relax\else
  \message{[already loaded]}\endinput\fi
\expandafter\edef\csname debug.def\endcsname{%
  \catcode`\noexpand\@\the\catcode`\@\edef\noexpand\styname{\styname}
  \def\expandafter\noexpand\csname debug.def\endcsname{\stydate}}
\catcode`\@=11 {\edef\temp{\the\everyjob\W@{\styname: <\stydate>}}
\global\everyjob\expandafter{\temp}}

\def\n@te#1#2{\leavevmode\vadjust{%
 {\setbox\z@\hbox to\z@{\strut\eightpoint\let\quotei\filename#1}%
  \setbox\z@\hbox{\raise\dp\strutbox\box\z@}\ht\z@=\z@\dp\z@=\z@%
  #2\box\z@}}}
\def\leftnote#1{\n@te{\hss#1\quad}{}}
\def\rightnote#1{\n@te{\quad\kern-\leftskip#1\hss}{\moveright\hsize}}
\def\?{\FN@\qumark}
\def\qumark{\ifx\next"\DN@"##1"{\leftnote{\rm##1}}\else
 \DN@{\leftnote{\rm??}}\fi{\rm??}\next@}
\def\filename#1{\hbox{\tt #1}}
\def\mnote@@#1{\rightnote{\vtop{%
 \ifcat\noexpand"\noexpand~\def"##1"{\filename{##1}}\fi
 \hsize2.0in \baselineskip7\p@\parindent\z@
 \tolerance\@M\spaceskip2.6\p@ plus10\p@ minus.9\p@\rm#1}}}
\def\mnote#1{\@bsphack\mnote@{#1}\@esphack}

\def\nonotes{\let\mnote@\eat@}
\def\printnotes{\let\mnote@\mnote@@}
\printnotes

\def\PrintLabels{%
 \gdef\printlabel@##1##2{\ifvmode\else\leftnote{\eighttt##2}\fi}}
\def\NoLabels{\global\let\printlabel@\eat@ii}
\NoLabels \@@addparm\everylabel@{\printlabel@#1}

\def\PrintFiles{\gdef\outputmark@{\line{\hfill\smash{\raise1cm\vbox{%
  \hbox to\z@{\kern1cm\tenrm\the\month/\the\day/\the\year\hss}%
  \hbox to\z@{\kern1cm\tentt\jobname.tex\hss}%
  \ifx\filecomment\undefined\else\hbox to\z@{\kern1cm\tenrm\filecomment\hss}\fi}}}}}
\def\NoFiles{\global\let\outputmark@\empty}
\def\PageMark{\gdef\outputmark@}
\let\outputmark@\empty
\@@addto\pagetop@\outputmark@

\def\tracerefs{\def\wrn@reference##1##2{\W@X{\string\reference{##1}}}}
\def\tracecites{\printwarning\nocite
  \def\wrn@citation##1##2{\@Nxdef\cite@-##1\@{##1}\W@X{\string\citation{##1}}}}
\def\wrn@@nocite#1#2{\ifRefs\wrn@@@nocite{#2}\fi}
\def\wrn@@@nocite#1{\if\@undefined\cite@-#1\@\warning@{citation `#1' [\@lastmark] not used}\fi}
\let\wrn@nocite\eat@ii
\let\reference\eat@
\let\citation\eat@
\@@addparm\everylabel@{\wrn@nocite#1}

\@@addto\disablepreambule@cs{%
    \/\PrintFiles\empty
    \/\NoFiles\empty
    \/\PageMark\empty
}

\csname debug.def\endcsname

\def\stydate{November 27, 1997}
\immediate\write16{This is DEGT.DEF by A.Degtyarev <\stydate>}
{\edef\temp{\the\everyjob\immediate\write16{DEGT.DEF: <\stydate>}}
\global\everyjob\expandafter{\temp}}

\chardef\tempcat\catcode`\@\catcode`\@=11

\let\ge\geqslant
\let\le\leqslant
\def\C{{\Bbb C}}
\def\R{{\Bbb R}}
\def\Z{{\Bbb Z}}
\def\Q{{\Bbb Q}}

\let\ZZ\Z

\def\Cp#1{\C{\operator@font p}^{#1}}
\def\Rp#1{\R{\operator@font p}^{#1}}
\def\Hom{\qopname@{Hom}}
\def\Ext{\qopname@{Ext}}
\def\Tors{\qopname@{Tors}}

\def\Im{\qopname@{Im}}          
\def\Re{\qopname@{Re}}     %
\def\Ker{\qopname@{Ker}}
\def\Coker{\qopname@{Coker}}
\def\Int{\qopname@{Int}}
\def\Cl{\qopname@{Cl}}
\def\Fr{\qopname@{Fr}}
\def\Fix{\qopname@{Fix}}
\def\tr{\qopname@{tr}}
\def\inj{\qopname@{in}}
\def\id{\qopname@{id}}
\def\pr{\qopname@{pr}}
\def\rel{\qopname@{rel}}
\def\pt{{\operator@font{pt}}}
\def\const{{\operator@font{const}}}
\def\codim{\qopname@{codim}}
\def\cdim{\qopname@{dim_{\C}}}
\def\rdim{\qopname@{dim_{\R}}}
\def\conj{\qopname@{conj}}
\def\rank{\qopname@{rk}}
\def\sign{\qopname@{sign}}
\def\gcd{\qopname@{g.c.d.}}
\let\sminus\smallsetminus
\def\set<#1|#2>{\bigl\{#1\bigm|#2\bigr\}}


\def\preprint#1{\hrule height0pt depth0pt\kern-24pt%
  \hbox to\hsize{#1}\kern24pt}
\def\today{\ifcase\month\or January\or February\or March\or
  April\or May\or June\or July\or August\or September\or October\or
  November\or December\fi \space\number\day, \number\year}

\def\n@te#1#2{\leavevmode\vadjust{%
 {\setbox\z@\hbox to\z@{\strut\eightpoint#1}%
  \setbox\z@\hbox{\raise\dp\strutbox\box\z@}\ht\z@=\z@\dp\z@=\z@%
  #2\box\z@}}}
\def\leftnote#1{\n@te{\hss#1\quad}{}}
\def\rightnote#1{\n@te{\quad\kern-\leftskip#1\hss}{\moveright\hsize}}
\def\?{\FN@\qumark}
\def\qumark{\ifx\next"\DN@"##1"{\leftnote{\rm##1}}\else
 \DN@{\leftnote{\rm??}}\fi{\rm??}\next@}

\def\centerpage{\dimen@=6.5truein \advance\dimen@-\hsize\hoffset.5\dimen@}
\ifnum\mag>1000 \centerpage\fi

\def\nologo{\let\logo@\relax}

\expandafter\ifx\csname eat@\endcsname\relax\def\eat@#1{}\fi
\expandafter\ifx\csname operator@font\endcsname\relax
 \def\operator@font{\roman}\fi
\expandafter\ifx\csname eightpoint\endcsname\relax
 \let\eightpoint\small\fi

\catcode`\@\tempcat\let\tempcat\undefined

{\catcode`\@=11 \gdef\proclaimfont@{\sl}}

\PrintLabels \NoLabels 

\let\GS\Sigma
\let\<\langle
\let\>\rangle

\def\CS{\Cal S}

\def\Cp#1{\Bbb P^{#1}}
\def\Rp#1{\Cp{#1}_\R}
\def\ZZ{\Z_2}
\def\+{\mathbin{\scriptstyle\sqcup}}
\def\ls|#1|{\mathopen|#1\mathclose|}
\def\Spin{\operatorname{Spin}}
\def\sbp{_{\scriptscriptstyle+}}

\newtheorem{}\claim\thm\endproclaim


\topmatter

\author
Alex Degtyarev and Viatcheslav Kharlamov
\endauthor

\title
Real Rational Surfaces Are Quasi-Simple
\endtitle

\address
Bilkent University, Turkey
\endaddress

\email
degt\,\@\,fen.bilkent.edu.tr
\endemail

\address
Institut de Recherche Math\'ematique Avanc\'ee,\newline\indent Universit\'e
Louis Pasteur et CNRS, France
\endaddress

\email
kharlam\,\@\,math.u-strasbg.fr
\endemail

\abstract We show that real rational (over~$\C$) surfaces are
quasi-simple, i.e., that such a surface is determined up to deformation
in the class of real surfaces by the topological type of its real
structure.
\endabstract



\endtopmatter

\document

\section{Introduction}\label{s-intro}

The list of all minimal models of complex rational surfaces is exhausted
by the projective plane~$\Cp2$ and the geometrically ruled surfaces
$\GS_a\to\Cp1$, $a\in\Z_+$. As is well known, this implies that
diffeomorphic rational surfaces (or, equivalently, rational surfaces with
isomorphic cohomology rings) are deformation equivalent. In the present
paper we treat the corresponding deformation problem for surfaces defined
over~$\R$ and rational over~$\C$ (which we call \emph{real rational
surfaces}, see the definition in Section~\ref{s-main}). We show that two
real rational surfaces are deformation equivalent over~$\R$ (i.e., in the
class of real surfaces) if and only if their real structures are
diffeomorphic. The prove is based on the classification of $\R$-minimal
real rational surfaces going back to A.~Comessatti~\cite{Co-rat}.

Note
that
the same phenomenon, which we call quasi-simplicity of the class of
surfaces, is observed for real Abelian surfaces (essentially due to
A.~Comessatti~\cite{Co-abel}), for real hyperelliptic surfaces
(F.~Catan\-ese and P.~Frediani~\cite{CF}), for real $K3$-surfaces
(essentially due to V.~Nikulin~\cite{NIK}), and for real Enriques
surfaces (A.~Degtyarev and V.~Kharlamov; the quasi-simplicity statement
was announced in~\cite{DK}, and the complete list of deformation classes
of real Enriques surfaces was obtained in collaboration with I.~Itenberg
in~\cite{Book}). It is natural to expect that such a simple behaviour
would no longer take place for more complicated surfaces, like those of
general type. Note, however, that the so called
\emph{Bogomolov-Miayoka-Yau surfaces}, which are of general type, are
also quasi-simple, see \cite{KK}. It is also worth mentioning that the
first non trivial quasi-simplicity result, concerning real cubic surfaces
in~$\Cp3$, was discovered by F.~Klein and
L.~Schl\"afli (see, for example, our survey~\cite{survey}). To our
surprise, apart from this result, we could not find any trace of
deformation classification of real rational surfaces in the literature.
(The topology of the real part of a real rational surface is discussed
in~\cite{Si}.)

It is worth
mentioning that
in all
cases above, including the case of rational surfaces, the deformation
class is determined, in fact,
by the {\bf topological} type of the real structure.
In particular, as a consequence, the topological type of the real
structure determines it up to diffeomorphism.

The quasi-simplicity of real rational surfaces was obtained during our
work on the deformation classification of real Enriques
surfaces. One particular case, concerning the
$\R$-minimal rational surfaces, is mentioned in~\cite{Segovia}. Another
one, concerning the Del Pezzo surfaces,
is contained in~\cite{Book}; this case is used explicitly in the
classification of real Enriques surfaces. Here we treat the remaining
cases and complete the proof.

\remark{Remark}
In spite of its quasi-simplicity (or due to some uncertainty principle?)
the class of real rational surfaces seems to be one of the most difficult
for study in what concerns other geometric properties. For example, it is
unknown whether the number of isomorphism classes of real structures on a
given rational surface is finite, although the corresponding finiteness
result is known to hold for all surfaces of Kodaira dimension $\ge1$ and
for minimal Abelian, hyperelliptic, $K3$-, and Enriques surfaces,
see~\cite{Book}.
\endremark

\subsection*{Acknowledgements}
We are grateful to J.~Koll\'ar for his helpful remarks and for his interest
in the quasi-simplicity phenomenon that gave us an additional stimulus
for publishing this proof. This paper was written during the first
author's stay at \emph{Universit\'e Louis Pasteur}, Strasbourg, supported
by the European Doctoral College and by a joint {\it CNRS\/}--{\it
T\"UB\accent95ITAK\/} exchange program.

\section{Main Result}\label{s-main}

\subsection{Basic definitions}
 From now on by a \emph{complex surface} we mean a compact complex
analytic variety of complex dimension~$2$. A \emph{real structure} on a
complex surface~$X$ is an anti-holomorphic involution $X\to X$. A surface
supplied with a real structure is called a \emph{real surface}; the real
structure is usually denoted by~$\conj$. The fixed point set
$\Fix\conj=X_\R$ is called the \emph{real part} of~$X$. Topologically,
$X_\R$ is a closed $2$-manifolds. To describe its homeomorphism type, we
will list its connected components in an expression of the form
$nS\+n_1S_1\+\dots\+m_1V_1\+\ldots$, where $S=S_0=S^2$ stands for the
$2$-sphere, $S_p=\#_p(S^1\times S^1)$, for the orientable surface of
genus~$g$, and $V_q=\#_q\Rp2$, for the nonorientable surface of
genus~$q$.

\remark{Remark}
Of course, if $X$ is an algebraic surface defined over~$\R$, then $\conj$
is the Galois involution and $X_\R$ is the set of $\R$-rational closed
points of~$X$. However, although all surfaces treated in this paper are
in fact algebraic, in general we prefer the more geometric complex
analytic approach.
\endremark

By a \emph{real rational surface} we mean a real surface birationally
equivalent to~$\Cp2$ {\bf over~$\C$}. Note that this terminology differs
from the commonly accepted one: `most' real rational surfaces in our
sense are not in fact rational over~$\R$; Theorem~\ref{Comessatti} below
gives a number of examples.

The meaning of many common notions and constructions depends on whether
they are considered over~$\C$ or over~$\R$. (The rationality of a real
surface is a good example.) The $\conj$-equivariant category is usually
designated with the word `real'; for example, a \emph{real curve} on a
real surface is a $\conj$-invariant analytic curve. However, sometimes
this approach is unacceptable; in these cases we will use the
prefix~$\R$. Thus, we will speak about \emph{$\R$-minimal} surfaces
(instead of `really minimal real surfaces'): a real surface~$X$ is called
$\R$-minimal if any $\conj$-invariant degree~$1$ holomorphic map $X\to Y$
to another real surface~$Y$ is a biholomorphism. Clearly, $X$ is
$\R$-minimal if and only if it is free of real $(-1)$-curves and of pairs
of disjoint $(-1)$-curves transposed by~$\conj$.

The following is the key notion in the topological study of
analytic varieties. The definition of deformation equivalence below
applies to both complex and real categories; in the former case the
prefix `real' should be omitted.

\definition
A \emph{deformation} of complex surfaces is a proper analytic submersion
$p\:Z\to D^2$, where $Z$ is a $3$-dimensional analytic variety and
$D^2\subset\C$ a disk. If $Z$ is real and $p$ is equivariant, the
deformation is called real. Two (real) surfaces~$X'$ and~$X''$ are called
\emph{deformation equivalent} if they can be connected by a chain
$X'=X_0$, \dots, $X_k=X''$ so that $X_i$ and $X_{i-1}$ are isomorphic
to (real) fibers of a (real) deformation.
\enddefinition

Unless stated otherwise, the deformation equivalence is understood in the
same category as the surfaces considered, i.e., real for real surfaces
and complex otherwise.

\definition\label{qs.def}
A real surface~$S$ is called \emph{quasi-simple} if it is deformation
equivalent to any other real surface~$S'$ such that \therosteritem1 $S'$
is deformation equivalent to~$S$ as a complex surface, and
\therosteritem2 the real structure of~$S'$ is diffeomorphic to that
of~$S$.
\enddefinition

Now we are ready to state our main result.

\theorem[\subsection{Main theorem}]\label{main}
Real rational surfaces are quasi-simple.
\endtheorem

In Section~\ref{s-nonminimal} below we reduce Theorem~\ref{main} to the
following statement, which is also of an independent interest.

\theorem[\subsection{Minimal surfaces}]\label{main-min}
Any two $\R$-minimal real rational surfaces with diffeomorphic real
structures are deformation equivalent in the class of $\R$-minimal real
rational surfaces.
\endtheorem

\remark{Remark}
In view of the modern advances in topology of $4$-manifolds one should
carefully distinguish between the smooth and topological categories, and
the former seems to us more appropriate in definitions like~\ref{qs.def}.
Note however, that in all known examples (see Introduction) the
quasi-simplicity holds in a stronger sense: homeomorphic real structures
(within a fixed complex deformation class) are deformation equivalent.
The same remark applies to rational surfaces: everywhere in this paper
one can replace the word `diffeomorphic' with `homeomorphic'
(and even `homotopy equivalent', as all the invariants
can be expressed in terms of homology, see Theorem~\ref{invariants}).
\endremark

\remark{Remark}
Alternatively, Theorem~\ref{main} can be strengthened in a different
direction.
Indeed, due to R.~Friedman and
Z.~Qin~\cite{Friedman}, a surface diffeomorphic to a rational one is
rational. Thus, one can state that, {\proclaimfont given a real rational
surface~$X$, any real surface~$X'$ diffeomorphic to~$X$ \rom(as a real
surface\rom) is deformation equivalent to~$X$.}
\endremark

\remark{Remark}
In~\ref{qs.def} we define the quasi-simplicity as a property of a real
surface or, equivalently, of a real deformation class. One could also
define as a property of a complex deformation class (i.e., as the
requirement that all the real surfaces that are $\C$-deformation
equivalent to a given one should be quasi-simple). In all examples known
so far the quasi-simplicity does hold for whole complex deformation
classes.
\endremark

\subsection{Plan of the proof}
The proof of Theorem~\ref{main-min} is based on the following
classification of $\R$-minimal real rational surfaces, due to
Comessatti~\cite{Co-rat} (see also~\cite{Si}, the survey~\cite{MT}, and
the lecture notes~\cite{Kollar}):

\theorem\label{Comessatti}
Each $\R$-minimal real rational surface~$X$ is one of the following\rom:
\roster
\item\local1
projective plane $\Cp2$ with its standard real structure\rom:
$X_\R=V_1$\rom;
\item\local2
quadric $\Cp1\times\Cp1$ with one of its four nonequivalent real
structures\rom: $X_\R=S_1$, $X_\R=S$, and two structures with
$X_\R=\varnothing$ \rom(see~\ref{rem.PxP} below\rom)\rom;
\item\local3
geometrically ruled surface $\GS_a$, $a\ge2$ \rom(see~\ref{rem.Sigma}
below\rom), with $X_\R=V_2$ and the standard real
structure, if $a$ is odd, and with $X_\R=S_1$ or~$\varnothing$ and one of
the two respective nonequivalent structures, if $a$ is even\rom;
\item\local4
relatively minimal real conic bundles over $\Cp1$ \rom(see
Section~\ref{s-bundles}\rom) with $2m\ge4$ singular fibers\rom:
$X_\R=mS$\rom;
\item\local5
unnodal Del Pezzo surfaces of degree $d=1$ or~$2$\rom: $X_\R=V_1\+4S$, if
$d=1$, and $X_\R=3S$ or $4S$, if $d=2$.
\endroster
Conversely, any rational surface from the above list is $\R$-minimal.
\qed
\endtheorem

\remark{Remark}
An $\R$-minimal unnodal Del Pezzo surface of degree~2 with $X_\R=3S$ can
also be represented as a conic bundle over~$\Cp1$ with six reducible
fibers. Note that \iref{Comessatti}5 does {\bf not} list all real
structures on (unnodal) Del Pezzo surfaces of degree~$1$ and~$2$; most of
them are not $\R$-minimal (see~\cite{Book, \S17.3}).
\endremark

To prove Theorem~\ref{main-min}, we first show that the statement holds
for each family listed in \ref{Comessatti} separately. The conic
bundles~\iref{Comessatti}4 are treated in Section~\ref{s-bundles}; the
other classes have already been considered elsewhere, see references in
Section~\ref{s-other}. Then it remains to show that whenever the
corresponding complex deformation classes overlap, the real surfaces are
deformation equivalent; this is done in~\ref{min-proof}.

\subsection{More notation}
Here we cite a few known results and introduce some notation used
throughout the paper.

\definition[]\label{rem.PxP}
Recall that the quadric $Y=\Cp1\times\Cp2$ has two nonisomorphic real
structures with $Y_\R=\varnothing$. They are both product structures:
$c_1\times c_1$ and $c_0\times c_1$, where $c_0$ is the ordinary complex
conjugation on~$\Cp1$ (with $\Rp1=S^1$) and $c_1$ is the quaternionic
real structure (with empty real point set). The quotient manifold
$Y/(c_1\times c_1)$ is not $\Spin$; the quotient $Y/(c_0\times c_1)$ is
$\Spin$. The real structure $c_1\times c_1$ is inherited from $\Cp3$
(with its standard real structure); we will call it \emph{standard}.

The two real structures with $Y_\R=\varnothing$ become equivalent after a
single blow-up at a pair of conjugate points (see, e.g.,~\cite{Book,
Lemma in~17.4.1}).
\enddefinition

\definition[]\label{rem.Sigma}
We denote by~$\GS_a$, $a\ge0$, the geometrically ruled rational surface
(i.e., relatively minimal conic bundle over~$\Cp1$, see
Section~\ref{s-bundles}) that has a section of square~$(-a)$. (Unless
$a=0$, such a section is unique; it is called the \emph{exceptional
section}.) All these surfaces except~$\GS_1$ are minimal. (Clearly,
$\GS_0=\Cp1\times\Cp1$ and $\GS_1$ is the plane $\Cp2$ blown up at one
point.) Recall that, up to isomorphism, $\GS_a$, $a>0$, has one real
structure if $a$~is odd (and then the real part is~$V_2$), and it has two
real structures if $a$~is even (the real part being either~$S_1$
or~$\varnothing$).

We will denote by~$\ls|e_0|$ and~$\ls|e_\infty|$ the classes of the
exceptional section and of generic section, respectively, so that
$e_0^2=-a$, $e_\infty^2=a$, and $e_0\cdot e_\infty=0$. (If $a=0$, then
$e_0=e_\infty$.) The class of the fiber will be denoted by~$\ls|l|$; one
has $l^2=0$ and $l\cdot e_0=l\cdot e_\infty=1$. As is well know, the
semigroup of classes of effective divisors on~$\GS_a$ is generated
by~$\ls|e_0|$ and~$\ls|l|$ and one has $\ls|e_\infty|=\ls|e_0+al|$.
\enddefinition

\definition[]
Given a topological space~$X$, we denote by $b_i(X)=\rank H_i(X)$ and
$\beta_i(X)=\dim H_i(X;\ZZ)$ its $\Z$- and $\ZZ$-Betti numbers,
respectively. The
corresponding total Betti numbers are denoted by
$b_*(X)=\sum_{i\ge0}b_i(X)$ and $\beta_*(X)=\sum_{i\ge0}\beta_i(X)$.
Recall the following Smith inequality: {\proclaimfont for any real
variety~$X$ one has}
$$
\beta_*(X_\R)\le\beta_*(X)\quad\text{and}\quad
\beta_*(X_\R)=\beta_*(X)\bmod2.
$$
Thus, one has $\beta_*(X_\R)=\beta_*(X)-2d$ for some nonnegative
integer~$d$; in this case $X$ is called an \emph{$(M-d)$-variety}.

If $X$ is a curve, the Smith inequality reduces to the so called Harnack
inequality: {\proclaimfont the number of connected components of the real
part of a nonsingular irreducible real curve~$X$ of genus~$g$ does not
exceed $g+1$.} Thus, $M$-curves of genus~$g$ are those with $(g+1)$
components in the real part.
\enddefinition

\definition[]
A \emph{cyclic order} on a finite set~$\CS$ is an embedding $\CS\to S^1$
considered up to auto-homeomorphism of the circle~$S^1$. In other words,
we consider cyclic orders up to reflections. Alternatively, a cyclic
order can be defined as a symmetric binary relation `to be neighbors'
such that each element of~$\CS$ has exactly two neighbors.
\enddefinition

\section{Conic Bundles}\label{s-bundles}

Recall that a \emph{conic bundle} is a regular map $p\:X\to C$, where
$X$~is a compact algebraic surface and $C$~is a curve, such that a
generic fiber of~$p$ is an irreducible rational curve. If both~$X$
and~$C$ are real and $p$ is equivariant, the conic bundle is called real.

As is well known, over~$\C$ any conic bundle blows down to the ruling of
a geometrically ruled surface. Over~$\R$, a real conic bundle is
relatively minimal if and only if all its singular fibers are real and
each singular fiber consists of a pair of conjugate $(-1)$-curves. As it
follows from~\ref{Comessatti}, a real conic bundle is relatively minimal
over~$\R$ if and only if it is $\R$-minimal.

\subsection{Birational classification}
An \emph{elementary transformation} of a (complex) conic bundle $X\to C$
at a point $P$ in a nondegenerate fiber~$F$ of~$X$ is the blow-up of~$P$
followed by the blow-down of~$F$. If $X\to C$ is real, a \emph{real
elementary transformation} of~$X$ is either an elementary transformation
at a real point of~$X$ or a pair of elementary transformations at two
conjugate points $P_1,P_2\in X$ that belong to two distinct conjugate
fibers of $X\to C$.

Next two statements are found in~\cite{Si} and~\cite{MT}.

\theorem\label{bundle-maps}
Any \rom(real\rom) fiberwise birational map between two relatively
minimal \rom(real\rom) conic bundles decomposes into a product of
\rom(real\rom) elementary transformations.
\qed
\endtheorem

\theorem\label{bundle-class}
Any conic bundle $X\to\Cp1$ is fiberwise birationally equivalent to a
conic bundle $Y\to\Cp1$ which is a branched double covering of a
geometrically ruled surface. Furthermore, $Y$ can be chosen so that in
appropriate affine coordinates $(x,y,z)$ it is given by the equation
$x^2+y^2=p(z)$, where $p$ is a polynomial, and the projection $Y\to\Cp1$
is the restriction to~$Y$ of the map $(x,y,z)\mapsto z$.

If $X\to\Cp1$ is real and $X_\R\ne\varnothing$, the above statement holds
in the equivariant setting, i.e., the bundle $Y\to\Cp1$, the coordinates
and the equation $x^2+y^2=p(z)$, and the birational equivalence can be
chosen real.
\qed
\endtheorem

The following is straightforward:

\lemma\label{Sigma}
Let a real conic bundle~$Y$ be the double covering of a geometrically
ruled surface $\GS_a$ branched over a nonsingular curve~$B$. Assume that
$Y$ has $2m\ge4$ singular fibers. Then one has $\ls|a|\le m$ and
$a=m\bmod2$ and the branch curve belongs to the linear system
$\ls|2e_\infty+(m-a)l|$. Furthermore, $Y$ is $\R$-minimal if and only if
$B$ is an $M$-curve.
\qed
\endlemma

\theorem[\subsection{Double ruled surfaces}]\label{th.double}
Two $\R$-minimal real rational conic bundles $Y_1$, $Y_2$ which are
double branched covering of geometrically ruled surfaces are deformation
equivalent if and only if they have the same number of singular fibers.
\endtheorem

\proof
We concentrate on the `if' part of the statement; the `only if' part is
obvious.

If $Y_1$ and $Y_2$ cover the same geometrically ruled surface~$\GS_a$,
then the statement follows from the classification of real curves in the
linear system $\ls|2e_\infty+(m-a)l|$ (see, e.g.,~\cite{survey}
or~\cite{Book, A2.2.8}): any two such $M$-curves are rigidly isotopic.

Thus, in view of~\ref{Sigma}, it remains to deform an $M$-curve in
$\ls|2e_\infty+(m-a)l|$ on~$\GS_a$ to an $M$-curve in
$\ls|2e_\infty+(m-a-2)l|$ on~$\GS_{a+2}$, $a\le m-2$. This can be done as
follows (cf.~\cite{Welsch}). Consider a generic $M$-curve
$C\in\ls|2e_\infty+(m-a-1)l|$ on the surface~$\GS_{a+1}$. Pick a real
point~$P_0$ of intersection of~$C$ and the exceptional section~$E_0$
of~$\GS_{a+1}$. (Such a point exists since $m-a-1>0$ is odd.) The desired
deformation is parametrized by the real points~$P_t$ of the germ of~$C$
at~$P_0$. It consists of the surfaces obtained from~$\GS_{a+1}$ by the
elementary transformations at~$P_t$ (and the corresponding transforms
of~$C$). If $t=0$, the resulting surface is $\GS_{a+2}$, otherwise it
is~$\GS_a$.
\endproof

\theorem[\subsection{The general case}]\label{th.bundle}
Two $\R$-minimal real rational conic bundles $Y_1$, $Y_2$ are deformation
equivalent if and only if they have the same number of singular fibers.
\endtheorem

The statement is an immediate consequence of~\ref{th.double} and the
following lemma.

\lemma
Any $\R$-minimal real rational conic bundle is deformation equivalent to
a conic bundle obtained by a double covering of a geometrically ruled
surface.
\endlemma

\proof
In view of~\ref{bundle-class} and~\ref{bundle-maps}, any real rational
conic bundle is obtained from a double covering by a sequence of real
elementary transformations. Up to deformation (moving the blow-up
centers), one can assume that the transformations are independent, and it
suffices to prove that a real elementary transformation applied to a
double covering results in a surface deformation equivalent to a double
covering.

Two complex conjugate blow-up centers can be merged to a single double
real center (recall that the real part of the surface is nonempty), and
then taken apart to two real centers. Thus, all the blow-up centers can
be assumed real. In view of~\ref{Sigma}, each component of the real part
of the surface intersects the branch curve of the covering; moving the
blow-up centers to the branch curve, one obtains a sequence of real
elementary transformations respecting the deck translation of the
covering. Hence, the resulting surface possesses an involution, i.e., it
is a double covering.
\endproof

\section{Other Minimal Surfaces}\label{s-other}

In this section we complete the proof of~\ref{main-min}. In fact, we
prove the following slightly stronger statement.

\theorem[\subsection{Theorem}]\label{main-min'}
With two exceptions, the real part $X_\R$ of an $\R$-minimal real
rational surface~$X$ determines~$X$ up to deformation in the class of
$\R$-minimal real rational surfaces. The exceptions are\rom:
\roster
\item
the two real structures on $X=\Cp1\times\Cp1$ with $X_\R=\varnothing$,
which differ by whether $X/\!\conj$ is $\Spin$ or not\rom;
\item
surfaces with $X_\R=4S$, see~\iref{Comessatti}4 and~\ditto5, which differ
by the topology of~$X$ \rom(the conic bundles in~\iref{Comessatti}4 are
$(M-2)$-surfaces, whereas the Del Pezzo surfaces of degree~$2$
in~\iref{Comessatti}5 are $(M-1)$-surfaces\rom).
\endroster
\endtheorem

\subsection{Proof of~\ref{main-min'}}\label{min-proof}
First, note that each family listed in~\ref{Comessatti} consists of
deformation equivalent surfaces. Indeed, the classification of real
structures on~$\Cp2$ and geometrically ruled surfaces
(\iref{Comessatti}1--\ditto3) is well known and reduces mainly to college
linear algebra, the classification of $\R$-minimal real rational conic
bundles (\iref{Comessatti}4) is given in~\ref{th.bundle}, and the
deformation classification of real Del Pezzo surfaces
(\iref{Comessatti}5) can be found, e.g., in~\cite{Book, \S17.3}. Note
also that {\bf all} surfaces listed in~\ref{Comessatti} are $\R$-minimal;
hence, their classification up to real deformation coincides with the
classification up to deformation in the class of $\R$-minimal real
surfaces. Thus, to complete the proof it remains to show that, with the
exception of the two cases mentioned in the statement, the distinct
families of surfaces with homeomorphic real parts are, in fact,
deformation equivalent. Below we consider these families case by case.

\subsubsection*{Cases 1, 2\rom: $X=\GS_a$, $X_\R=S_1$ if $a$ is even or $X_\R=V_2$ if $a$ is odd}
A real deformation $\GS_a\mapsto\GS_{a+2}$ can be constructed as in the
proof of~\ref{th.double}, with the branch curve ignored, i.e., both the
surfaces occur in a family of real elementary transformations of the
surface~$\GS_{a+1}$ (with nonempty real part).

\subsubsection*{Case 3\rom: $X=\GS_a$ with $a$ even, $X_\R=\varnothing$}
(We do not exclude the case when $X$ is an empty quadric in~$\Cp3$ with
its standard real structure.) A deformation $\GS_a\mapsto\GS_{a+4}$ can
be constructed as above, as a family of elementary transformations
of~$\GS_{a+2}$ with complex conjugate blow-up centers. For a deformation
$\GS_0\mapsto\GS_2$ one can take the family of surfaces given in
$\Cp1\times\Cp1\times\Cp1$ by the equations
$$
v_1^2y_1z_0 = v_0^2y_0z_1 + tv_0v_1(y_0z_0+y_1z_1),
$$
where $[v_0:v_1]$, $[y_0:y_1]$, $[z_0:z_1]$ are homogeneous coordinates
in the factors and the real structure is
$([v_0:v_1],[y_0:y_1],[z_0:z_1])\mapsto([-\bar v_1:\bar v_0],[-\bar
y_1:\bar y_0],[-\bar z_1:\bar z_0])$.
\qed

\section{Canonical Minimal Models}

In order to prove Theorem~\ref{main}, we consider a real rational
surface~$X$ and show that, up to deformation, $X$ has a certain canonical
minimal model whose topology and, hence, deformation type is determined
by the topology of the real structure of~$X$. Then, given two such
surfaces $X_1$, $X_2$ with diffeomorphic real structures, one can assume
that they are both obtained by a sequence of blow-ups from the same
$\R$-minimal surface~$Y$, and it would remain to show that the sets of
blow-up centers in~$Y$ are also determined by the topology of~$X_i$ up to
equivariant isotopy in~$Y$. This is done in Section~\ref{s-nonminimal}.

\subsection{Essential invariant classes}
Given a real rational surface $(X,\conj)$, denote by $L\sbp(X)$ the
eigenlattice
$$
L\sbp(X)=\Ker\bigl[(1-\conj_*)\:H_2(X)\to H_2(X)\bigr]
$$
supplied with the induced intersection form.

\lemma\label{negative}
The lattice $L\sbp(X)$ is negative definite and even. If another real
rational surface~$X'$ is obtained from~$X$ by a sequence of $r$ blow-ups
at real points and $q$ blow-ups at pairs of complex conjugate points,
then $L\sbp(X')=L\sbp(X)\oplus q\<-2\>$.
\endlemma

\proof
One has $\sigma\sbp H_2(X)=1$. By the usual averaging arguments one can
find a $\conj$-invariant K\"ahler metric on~$X$. The Poincar\'e dual
$\gamma\in H_2(X;\R)$ of the fundamental class of such a metric belongs
to the $(-1)$-eigenspace $\Ker(1+\conj_*)$. Hence, the latter contains a
vector with positive square, and its orthogonal complement
$L\sbp(X)\otimes\R$ is negative definite.

The fact that $L\sbp(X)$ is even is a well known common property of real
surfaces. It follows, for example, from the fact that the Stiefel-Whitney
class $w_2(X)$ (which coincides with the $\operatorname{mod}2$-reduction
of the first Chern class~$c_1(X)$\,) can be realized by an algebraic
cycle which is defined over~$\R$ and thus belongs to the
$(-1)$-eigenlattice of~$\conj_*$.

The other assertions of the lemma are straightforward.
\endproof

As is well known, a definite integral lattice admits a unique (up to
reordering) decomposition into irreducible orthogonal summands. Denote by
$E\sbp(X)$ the \emph{essential part} of $L\sbp(X)$ obtained by removing
from $L\sbp(X)$ the summands isomorphic to~$\<-2\>$. Due to
Lemma~\ref{negative} this lattice is well defined and depends only on the
topology of $(X,\conj)$. The following is also an immediate consequence
of Lemma~\ref{negative}.

\corollary\label{e-part}
The lattice $E\sbp(X)$ is a birational invariant of a real rational
surface~$X$, i.e., any real birational transformation $X_1\dasharrow X_2$
of real rational surfaces induces an isomorphism $E\sbp(X_1)=E\sbp(X_2)$.
\qed
\endcorollary

\theorem[\subsection{Canonical minimal models}]\label{min.model}
A real rational surface~$X$ is $\R$-deformation equivalent to a surface
obtained by a sequence of real blow-ups from an $\R$-minimal real
rational surface~$Y$ described below.
\roster
\item\local{b0}
If $b_0(X_\R)=0$ and $b_*(X)>4$, then $Y=\Cp1\times\Cp1$ with
the standard real structure with $Y_\R=\varnothing$
\rom(see~\ref{rem.PxP} and remark below\rom).
\item\local{b1or}
If $b_0(X_\R)=1$ and $X_\R$ is orientable, then $Y=\Cp1\times\Cp1$ with
$Y_\R=S$ or $S_1$.
\item\local{b1non}
If $b_0(X_\R)=1$ and $X_\R$ is nonorientable, then $Y=\Cp2$.
\item\local{b23}
If $b_0(X_\R)=2$ or~$3$, then $Y$ is an $\R$-minimal conic bundle.
\item\local{b4}
If $b_0(X_\R)=4$, then $Y$ is a conic bundle or a Del Pezzo surface of
degree~$2$. The two cases are distinguished by $E\sbp(X)=D_8$ or~$E_7$,
respectively.
\item\local{b5}
If $b_0(X_\R)=5$, then $Y$ is a conic bundle or a Del Pezzo surface of
degree~$1$. The two cases are distinguished by $E\sbp(X)=D_{10}$
or~$E_8$, respectively.
\item\local{b6}
If $b_0(X_\R)\ge6$, then $Y$ is an $\R$-minimal conic bundle.
\endroster
Thus, the deformation class of~$X$ admits a canonical $\R$-minimal model
whose deformation type is determined by the topology of $(X,\conj)$.
\endtheorem

\remark{Remark}
In fact, in~\iref{min.model}{b0} one can choose and fix any of the two
real structures on~$Y$ with $Y_\R=\varnothing$, as they become equivalent
after a single blow-up at a pair of conjugate points.
\endremark

\remark{Remark}
If $b_0(X_\R)=0$ and $b_*(X)=4$ (cf.~\iref{min.model}{b0}), then $X$ is
$\R$-minimal, see~\ref{Comessatti}. Hence, in this case $X$ is also
deformation equivalent to $Y=\Cp1\times\Cp1$, but the real structure
on~$Y$ depends on~$X$.
\endremark

\proof
Since the number of connected components of the real part of a surface is
obviously preserved by a blow-up, the possible minimal models of~$X$ can
be found using Theorem~\ref{Comessatti}. These models can further be
deformed using Theorem~\ref{main-min}. This completes the proof in all
cases except~\loccit{b0} and~\loccit{b1non}. (In cases~\loccit{b4}
and~\loccit{b5} the lattice $E\sbp(X)$ is found by a direct calculation.
In view of Corollary~\ref{e-part} it suffices to calculate the lattice
for the minimal models. In fact, for the two Del Pezzo surfaces in
question one has $E_+=K^\perp$, where $K$ is the canonical class, and the
latter lattice is well known, see, e.g.,~\cite{Manin}
or~\cite{Demazure}.)

In case~\loccit{b0} Theorems~\ref{Comessatti} and~\ref{main-min} yield
$Y=\Cp1\times\Cp1$ with one of the two real structures with empty real
part. The two structures become isomorphic after a blow-up at a pair of
conjugate points (see, e.g.,~\cite{Book, Lemma in 17.3.1}); hence, if $X$
is not $\R$-minimal, one can choose any structure.

In case~\loccit{b1non} Theorems~\ref{Comessatti} and~\ref{main-min} yield
either $Y=\Cp2$, or $Y=\GS_{2k+1}$, or $Y=\Cp1\times\Cp1$ with
$Y_\R\ne\varnothing$. The ruled surfaces $\GS_{2k+1}$ are deformation
equivalent to~$\GS_1$ (cf.~\ref{min-proof}, Cases~1,~2), which blows down
to~$\Cp2$. If $Y=\Cp1\times\Cp1$, then the passage to~$X$ contains at
least one blow-up at a real point (recall that $X_\R$ is nonorientable);
after such a blow-up the surface also blows down to~$\Cp2$.
\endproof

\section{Proof of the Main Theorem}\label{s-nonminimal}

\subsection{Cyclic order of components}\label{s-cyclic}
Let $Y\to\Cp1$ be an $\R$-minimal real rational conic bundle, $Y_\R=mS$.
Then the projections of the connected components of~$Y_\R$ are disjoint
segments in $\Rp1=S^1$. This defines a canonical cyclic order on the set
of components of~$Y_\R$. From Proposition~\ref{cyclic} below it follows,
in particular, that this cyclic order does not depend on the conic bundle
structure of~$Y$.

\proposition\label{cyclic}
Let $X$ be a real rational surface obtained by a sequence of real
blow-ups from an $\R$-minimal real rational conic bundle $Y\to\Cp1$. Then
the cyclic order of the components of~$Y_\R$ is determined by the classes
realized by the components of~$X_\R$ in $H_2(X/\!\conj{};\ZZ)$. In
particular, the set of connected components of $X_\R$ also has a natural
cyclic order.
\endproposition

\proof
Recall that for any real surface $X$ the projection $X\to X/\!\conj$
induces an isomorphism $L\sbp(X)\otimes\Q=H_2(X/\!\conj{};\Q)$. Hence,
from Lemma~\ref{negative} it follows that $H_2(Y/\!\conj)$ is a negative
definite lattice and an easy calculation for the minimal models shows
that $H_2(Y/\!\conj)=2m\<-1\>$. (Alternatively, this statement can be
derived from Donaldson's theorem~\cite{Donald} or from the fact that
$Y/\!\conj$ decomposes into connected sum of several copies of
$\overline{\Cp2}$, see~\cite{Finashin}.) Thus, $H_2(Y/\!\conj)$ has a
canonical (up to reordering and multiplication by~$(-1)$\,) basis
$e_1,\dots,e_{2m}$ such that $e_i^2=-1$ and $e_ie_j=0$ for $i\ne j$. This
basis can easily be visualized. Let $F_1,\dots,F_m$ be the components
of~$Y_\R$ (numbered according to their cyclic order) and $f_i\in
H_2(Y/\!\conj)$ their classes. Then one can take for $e_{2i+1}$ the class
realized by the image in $Y/\!\conj$ of a vanishing cycle in~$Y$
corresponding to joining the components~$F_i$ and $F_{i+1}$. (To unify
the notation, we use the `cyclic' numbering of the components and
classes: $F_{i+km}=F_i$, $e_{i+2km}=e_i$, etc.) Let, further,
$g_1,\dots,g_{2m}\in H_2(Y/\!\conj)$ be the classes realized by the
images of the singular fibers of~$Y$, numbered so that $F_i$ is adjacent
to the fibers corresponding to $g_{2i}$ and $g_{2i+1}$. One has
$g_i^2=-2$ and $g_ig_j=0$ for $i\ne j$. It follows that these elements
generate the lattice over~$\Q$, and calculating the intersection numbers,
one obtains $e_{2i+1}=\frac12(g_{2i+2}-g_{2i+1})$. Hence, the other
$m$~generators have the form $e_{2i+2}=\frac12(g_{2i+2}-g_{2i+1})$ and,
under appropriate choice of the signs, one has
$f_i=e_{2i+2}-e_{2i+1}-e_{2i}-e_{2i-1}$.

Since $H_2(Y/\!\conj)$ is torsion free (as $Y_\R\ne\varnothing$), the
group $H_2(Y/\!\conj{};\ZZ)=H_2(Y/\!\conj)\otimes\ZZ$ also has a
canonical basis $\bar e_i=e_i\bmod2$, $i=1,\dots,2m$ and, denoting by
$\bar f_i=f_i\bmod2$ the class realized by~$F_i$, one has $\bar f_i=\bar
e_{2i+2}-\bar e_{2i+1}-\bar e_{2i}-\bar e_{2i-1}$. Hence, the cyclic
order is determined as follows: {\proclaimfont two components $F_i$,
$F_j$ are neighbors if and only if there is a standard basis
element~$\bar e$ such that $\bar e\bar f_i=\bar e\bar f_j=1$.}

It remains to note that the above rule extends literally to any
surface~$X$ obtained from~$Y$ by a sequence of blow-ups. Indeed, since
$H_2(X/\!\conj)=H_2(Y/\!\conj)\oplus q\<-1\>$ (the extra summands
corresponding to the pairs of complex conjugate exceptional divisors),
this lattice and hence the group
$H_2(X/\!\conj{};\ZZ)=H_2(X/\!\conj)\otimes\ZZ$ has a canonical basis,
which can be used to determine the cyclic order.
\endproof

\proposition\label{cyclic.constr}
Any cyclic permutation and/or reversing of the order of the components of
the real part of an $\R$-minimal real rational conic bundle can be
realized by a deformation and automorphism of the surface.
\endproposition

\proof
The statement can easily be proved by constructing explicit symmetric
surfaces. A cyclic permutation sending each component to one of its
neighbors can be realized in the family~$Y_t$ of surfaces given by the
affine equation $x^2+y^2=(\omega z)^{m}+(\omega z)^{-m}$, where
$\omega=\exp(2\pi it/m)$, $t\in[0,1]$ is the parameter of the
deformation, and the real structure is given by $(x,y,z)\mapsto(\bar
x,\bar y,\bar z^{-1})$. The order of the components is reversed by the
automorphism $(x,y,z)\mapsto(x,y,z^{-1})$ of~$Y_0$.
\endproof

Let now $Y$ be a real degree~$1$ Del Pezzo surface with $Y_\R=V_1\+4S$.
Then $Y$ can canonically be represented as the covering of a quadratic
cone branched over the vertex and a real $M$-curve $B\in\ls|3e_\infty|$
(see, e.g., \cite{Demazure}; since the map is defined by the
anti-bicanonical system, it is necessarily real),
the spherical components of~$Y_\R$ corresponding to the ovals of
the branch curve. Hence, the ruling of the cone determines a canonical
cyclic order on the set of the spherical components of~$Y_\R$. Next two
statements are proved similar to~\ref{cyclic} and~\ref{cyclic.constr}.

\proposition\label{DP.cyclic}
Let $X$ be a real rational surface obtained by a sequence of real
blow-ups from a real degree~$1$ Del Pezzo surface $Y$ with
$Y_\R=V_1\+4S$. Then the cyclic order of the spherical components
of~$Y_\R$ and the $V_1$-component of~$Y_\R$ are determined by the classes
realized by the components of~$X_\R$ in $H_2(X;\ZZ)$ and
$H_2(X/\!\conj{};\ZZ)$. In particular, $X_\R$ has a distinguished
component that blows down to~$V_1$, and the set of the other components
of $X_\R$ has a natural cyclic order.
\qed
\endproposition

\remark{Remark}
The rule for determining the cyclic order of the components is the same
as in the case of conic bundles. The component of~$X_\R$ that blows down
to~$V_1$ can be distinguished as the `essentially nonorientable' one: it
is the only component whose class in
$$
\Ker\bigl[(1-\conj_*)\:H_2(X;\ZZ)\to H_2(X;\ZZ)\bigr]
$$
does not belong to the image of $L\sbp(X)$.
\endremark

\proposition\label{DP.constr}
Any cyclic permutation and/or reversing of the order of the spherical
components of the real part of a real degree~$1$ Del Pezzo surface~$Y$
with $Y_\R=V_1\+4S$ can be realized by a deformation and automorphism of
the surface.
\qed
\endproposition

\subsection{Permutations of components}
Let us show that the connected components of the real part of all other
$\R$-minimal real rational surfaces can be permuted by a deformation.

\proposition\label{permutations}
Let $Y$ be either an $\R$-minimal real rational conic bundle with
$Y_\R=2S$ or~$3S$ or a real degree~$2$ Del Pezzo surface with $Y_\R=4S$.
Then any permutation of the components of~$Y_\R$ can be realized by a
deformation.
\endproposition

\proof
A surface as in the statement with $m=2$, $3$, or~$4$ connected
components can be constructed from a plane $M$-quartic~$B$ with $(4-m)$
isolated double points: one blows up the double points of~$B$, blows down
the line through them (in the case $m=2$), and constructs the double
covering of the resulting surface branched over the transform of~$B$. The
statement follows now from the well known fact that the ovals of an
$M$-quartic can be permuted by a rigid isotopy (which, in turn, can
easily be proved by constructing symmetric quartics and using the
connectedness of the group $\text{\sl PGL\,}(3,\R)$\,).
\endproof

\subsection{Proof of Theorem~\ref{main}}
Let $X_1$, $X_2$ be two real rational surfaces with diffeomorphic real
structures. In view of~\ref{min.model}, $X_1$ and~$X_2$ are deformation
equivalent to some surfaces~$X_1'$, $X_2'$ that blow down to the same
$\R$-minimal model~$Y$. Let $R_i\subset Y_\R$ and $I_i\subset Y\sminus
Y_\R$ be the corresponding sets of blow-up centers, $i=1,2$. (Up to
deformation one can assume that all blow-ups are independent.) Clearly,
the numbers $\#R_1=\#R_2$ and $\#I_1=\#I_2$ are determined by the
topology of the involutions, and since $Y\sminus Y_\R$ is connected, one
can assume that $I_1=I_2$.

Let $F_1$, \dots, $F_k$ be the components of~$Y_\R$.
From~\ref{cyclic}--\ref{DP.constr} and~\ref{permutations} it follows that
the components of one of the two copies of~$Y$ can be reordered by a
deformation and automorphism of the surface so that for each
$j=1,\dots,k$ one has $\#(R_1\cap F_j)=\#(R_2\cap F_j)$. Then one can
also isotope $R_1$ to~$R_2$, which results in a deformation of~$X_1'$
to~$X_2'$.
\qed

\subsection{Invariants}
In fact, we have proved that the deformation type of a real rational
surface is determined by certain homological information. More precisely,
extracting the invariants used in the proofs (see~\ref{min.model}
and~\ref{s-cyclic}), one arrives at the following statement.

\theorem\label{invariants}
Two real rational surfaces $X'$, $X''$ are deformation equivalent if and
only if there is a homeomorphism $X'_\R\to X''_\R$ preserving the
additional structures described below and one has $E\sbp(X')=E\sbp(X'')$
\rom(see~\ref{e-part}\rom; this condition holds automatically unless
$b_0(X'_\R)=b_0(X''_\R)=4$ or $5$\rom).

The additional structures to be preserved are\rom:
\roster
\item
a distinguished `essentially nonorientable' component of~$X_\R$, defined
in the case $E\sbp(X)=E_8$\rom: this is the only component whose class in
$H_2(X;\ZZ)$ does not belong to the image of $L\sbp(X)$\rom;
\item
a cyclic order on the set of the other components, defined whenever
$b_0(X)\ge4$ as follows\rom: two components $F_1$, $F_2$ are neighbors if
and only if there is a canonical basis element $\bar e\in
H_2(X/\!\conj;\ZZ)=H_2(X/\!\conj)\otimes\ZZ$ \rom(see the proof
of~\ref{cyclic}\rom) such that $\bar e[F_1]=\bar e[F_2]=1$.
\qed
\endroster
\endtheorem

\enddocument